\documentclass[11pt]{article}
\usepackage{amssymb,latexsym,amsmath}

\def\Z{\mathbb{Z}}
\def\N{\mathbb{N}}

\def\R{\mathbb{R}}
\def\C{\mathbb{C}}
\def\I{{\cal I}}
\def\P{{\cal P}}
\def\RR{{\cal R}}
\def\T{{\cal T}}
\def\s{s} 

\begin{document}
\setlength{\parindent}{0pt}
\setlength{\parskip}{0.4cm}

\newtheorem{theorem}{Theorem}[section]
\newtheorem{lemma}[theorem]{Lemma}
\newtheorem{proposition}[theorem]{Proposition}
\newtheorem{corollary}[theorem]{Corollary}
\newtheorem{definition}[theorem]{Definition}

\begin{center}
\large{\bf Explicit and efficient formulas for the lattice point count in rational polygons using Dedekind--Rademacher sums} 
\footnote{Appeared in \emph{Discrete and Computational Geometry} {\bf 27} (2002), 443-459. \\ 
          Parts of this work appeared in the first author's Ph.D. thesis.  
          The second author kindly acknowledges the support of NSA grant MSPR-OOY-196. 
          } 

\normalsize{\sc Matthias Beck and Sinai Robins} 

\end{center}

\footnotesize {\bf Abstract.} We give explicit, polynomial--time computable formulas for 
the number of integer points in any two--dimensional rational polygon. A rational polygon 
is one whose vertices have rational coordinates.  
We find that the basic building blocks of our formulas 
are \emph{Dedekind--Rademacher sums}, which are polynomial--time computable finite Fourier series. 
As a by--product we rederive a reciprocity law for 
these sums due to Gessel, which generalizes the reciprocity law for the classical Dedekind sums. 
In addition, our approach shows that Gessel's reciprocity law is a special case of the one for  
Dedekind--Rademacher sums, due to Rademacher.  


\vspace{5mm}

\small
\begin{quote} 
{\it The full beauty of the subject of generating functions emerges only from tuning in on both channels: the discrete and the continuous.} 
Herb Wilf \cite[p. vii]{wilfgen} 
\end{quote} 
\normalsize


\section{Introduction} 
We define a \emph{two--dimensional polytope} $ \P $ as a compact subset of $ \R^{ 2 } $ bounded by a simple, 
closed polygonal curve. $ \P $ is called a \emph{rational polytope} if all of its vertices have rational 
coordinates. We give explicit, polynomial--time computable (in the logarithm of the 
coordinates of the vertices) formulas for the number of 
integer points in any two--dimensional rational polytope and its integral dilations. 
We emphasize an expository flavor in this paper. In the current literature, there are either 
`nice' formulas that do not appear to be polynomial--time computable \cite{brion,diaz,kantor,puk,pomm}, or there are polynomial--time 
computable algorithms without `nice' formulas \cite{barvinok}. Asking for both seems to be asking for too much; 
but in $ \R^{ 2 } $ we show that we can have our cake and eat it, too. 

To fix notation, let $ \P^{\circ} $ be the interior of $ \P $, and $ \overline{\P} = \P $ be the 
closure of $ \P $. For $ t \in \N $, let $ L ( \P^{\circ}  , t ) = \# \left( t \P^{\circ}  \cap \Z^{ 2 }  \right) $ 
and $ L ( \overline{\P} , t ) = \# \left( t \overline{\P} \cap \Z^{ 2 }  \right) $ be 
the number of lattice points in the interior and closure, respectively, of the dilated polytope 
$ t \P = \{ ( tx , ty ) : (x,y) \in \P \} $. 
Ehrhart, who initiated the study of the lattice point count in dilated polytopes \cite{ehrhart}, 
proved that $ L ( {\cal P}^{\circ}  , t ) $ and $ L ( \overline{\cal P} , t ) $ are quasipolynomials 
in $t$. A \emph{quasipolynomial} is an expression of the form $ c_{n}(t) \ t^{n} + \dots + c_{1}(t) \ t + c_{0}(t) $, 
where $c_{0}, \dots , c_{n}$ are periodic functions in $t$. 

A natural first step is to fix a triangulation of $ \P $, which reduces 
our problem to counting integer points in rational \emph{triangles}. 
However, this procedure merits some remarks.  First, triangulation 
is not easy, but there has been remarkable progress recently, so 
that we can triangulate a convex polygon with $n$ vertices in roughly 
$n$ steps \cite{chazelle}. On the other hand, we are concerned with the 
efficiency of our formulas with respect to the coordinates of the 
vertices of $ \P $, and not the number of vertices. 

Another point of non--trivial significance is the number of lattice 
points on rational line segments, namely the boundary of our 
triangles. Although this is considerable easier than enumerating 
lattice points in 2--dimensional regions, it is still non--trivial 
and has only recently been completely solved \cite{br,tripati}. It is 
amusing that counting lattice points on line segments gave rise to links 
with the Frobenius coin--exchange problem and the number of representations 
of an integer by a linear form. 

After triangulating $\P$, we can further simplify the picture by embedding an arbitrary rational triangle in a rational rectangle: 

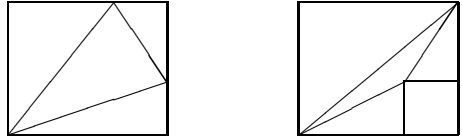
\begin{center} 
\begin{picture}(170,50)

\put(0,0){\framebox(60,50)}
\put(0,0){\line(4,5){40}} 
\put(40,50){\line(2,-3){20}}
\put(0,0){\line(3,1){60}} 

\put(110,0){\framebox(60,50)} 
\put(150,0){\framebox(20,20)} 
\put(110,0){\line(6,5){60}} 
\put(110,0){\line(2,1){40}} 
\put(150,20){\line(2,3){20}} 

\end{picture} 
\end{center} 

Since rectangles are easy to deal with, the problem reduces to finding a formula for a right--angled rational triangle. 
Such a formula is given in Section \ref{mattSec} using generating functions; this derivation is a refinement of a previously introduced method \cite{beck}. 
We find (Section \ref{sums}) that the basic building blocks of the lattice point count formulas for 
any two--dimensional rational polytope are the sawtooth function  
  \[ ((x)) := x - [x] - \frac{ 1 }{ 2 } \] 
and the {\it Dedekind--Rademacher sum} 
  \begin{equation}\label{s_t} \sigma (a,b,t) := \sum_{ k=0 }^{ b-1 } \left( \left( \frac{ ak+t }{ b } \right)  \right) \left( \left( \frac{ k }{ b }  \right)  \right) \ . \end{equation}  
Here $a$ and $b$ are integers, and $t$ is a real number.  
We use the name Dedekind--Rademacher sum in a somewhat lenient fashion; often $ ((x)) $ is  
defined to be 0 if $ x \in \Z $, also Rademacher's original definition is \cite{dieter,meyer,rademacher} 
  \[ \s ( a,b;x,y ) := \sum_{ k=0 }^{ b-1 } \left( \left( \frac{ a(k+y) }{ b } + x \right)  \right) \left( \left( \frac{ k+y }{ b }  \right)  \right) \ . \]  
Here $a$ and $b$ are integers, whereas $x$ and $y$ are real.  
However, it is clear that the different use of $(( \dots ))$ only results in a difference of  
the arithmetic sums in a trivial term. Also, $ \sigma $ and $ \s $ are strongly linked via  
  \[ \sigma (a,b,t) = \s \left( a, b; \frac t a , 0 \right) \]  
and  
  \[ \s ( a,b;x,y ) = \sigma ( a,b, ay+bx ) + \frac y b (( ya + xb )) \ . \]  
We chose to use $ \sigma $ rather than $ \s $ because of its natural appearance in our formulas.  
There exists a two--term reciprocity law for these sums \cite{knuth,rademacher}, which enables us 
to compute $ \sigma (a,b,t) $ in polynomial time, similar in spirit to the Euclidean 
algorithm. From this fact we conclude that our lattice point enumerator for $ \P $ is polynomial--time computable (Section \ref{consequences}). 

As a by--product of our formulas, we rederive in Section \ref{laws} two reciprocity laws as corollaries: the two--term law for the 
classical Dedekind sum (\cite{dedekind}, Chapter 2 of \cite{grosswald}), and a two--term law for generalized Dedekind 
sums due to Gessel \cite{gessel}. 
In fact, our approach shows that Gessel's reciprocity law is a special case of the reciprocity  
law for Dedekind--Rademacher sums, a theorem due to Rademacher \cite{rademacher}.  


\section{Generating functions}\label{mattSec} 
In \cite{beck}, the first author used the residue theorem to count lattice points in 
certain tetrahedra. Here we adjust and expand these methods to the rectangular triangles we reduced 
the discussion to in the introduction. Such a rectangular triangle $\T$ is given as a 
subset of $ \R^{2} $ consisting of all points $ (x,y) $ satisfying 
  \[ x \geq \frac{ a }{ d }  , \ y \geq \frac{ b }{ d } , \ e x + f y \leq r \]  
for some integers $ a,b,d,e,f,r $ with $ e a + f b  \leq r d $. 
Because the lattice point count is invariant under horizontal and vertical integer translation 
and under flipping about x- or y-axis, we may assume that $ a,b,d,e,f,r \geq 0 $ and $ a,b < d $. 
Let's further factor out the greatest common divisor $c$ of $e$ and $f$, so that $ e = cp $ and 
$ f = cq $, where $p$ and $q$ are relatively prime. Hence 
  \begin{equation}\label{triangle} \T = \left\{ (x,y) \in \R^{ 2 } : \ x \geq \frac{ a }{ d }  , \ y \geq \frac{ b }{ d } , \ cp x + cq y \leq r \right\} \ . \end{equation} 
To derive a formula for $ L \left( \overline{\T} , t \right) $ we interpret, 
similarly as in \cite{beck}, 
  \[ L \left( \overline{\T} , t \right) = \# \left\{ (m,n) \in \Z^{ 2 } : \ m \geq \frac{ ta }{ d }  , \ n \geq \frac{ tb }{ d } , \ cp m + cq n \leq tr \right\} \]
as the Taylor coefficient of $ z^{ t r  } $ of the function
\begin{eqnarray} &\mbox{}& \left( \sum_{ m \geq \left[ \frac{ t a - 1 }{ d  }  \right] + 1 } z^{ cp m }  \right) \left( \sum_{ n \geq \left[ \frac{ t b - 1 }{ d  }  \right] + 1 } z^{ cq n } \right) \left( \sum_{ k \geq 0 } z^{ k }  \right) = \frac{ z^{ \left( \left[ \frac{ t a - 1 }{ d  }  \right] + 1 \right) cp  }  }{ 1 - z^{ cp  }  } \frac{ z^{ \left( \left[ \frac{ t b - 1 }{ d  }  \right] + 1 \right) cq  }  }{ 1 - z^{ cq  }  } \frac{ 1 }{ 1 - z } \nonumber \\ 
                  &\mbox{}& \qquad = \frac{ z^{ u+v }  }{ \left( 1 - z^{ cp  } \right) \left( 1 - z^{ cq  } \right) ( 1 - z ) } \ , \label{function} \end{eqnarray} 
where we introduced, for ease of notation, 
  \begin{equation}\label{uandv} u := \left( \left[ \frac{ ta - 1 }{ d }  \right] + 1 \right) cp \qquad \mbox{ and } \qquad v := \left( \left[ \frac{ tb - 1 }{ d }  \right] + 1 \right) cq \ . \end{equation} 
We present two methods on how to extract the lattice point count from this generating function: 
{\it partial fractions} and the {\it residue theorem}. Both are inspired by works on generalized Dedekind sums, 
the first one by Gessel \cite{gessel}, the latter one by Zagier \cite{zagier}. In fact, both ways are completely 
equivalent, since our generating function is rational. However, to please both algebraically and analytically 
minded readers, we give two proofs of the following 
\begin{proposition}\label{closure} For the rectangular rational triangle $\T$ given by (\ref{triangle}), 
  \begin{eqnarray*} &\mbox{}& L \left( \overline{\T} , t \right) = \frac{ 1 }{ 2 c^{2} p q  } \left( t r - u - v \right)^{ 2 } + \frac{ 1 }{ 2 } \left( t r - u - v \right) \left( \frac{ 1 }{ cp  } + \frac{ 1 }{ cq  } + \frac{ 1 }{ c^{2} p q  }  \right) \\ 
                    &\mbox{}& \quad \qquad + \frac{ 1 }{ 4 } \left( 1 + \frac{ 1 }{ cp  } + \frac{ 1 }{ cq  }  \right) + \frac{ 1 }{ 12 } \left( \frac{ p  }{ q  } + \frac{ q  }{ p  } + \frac{ 1 }{ c^{2} p q  }  \right) \\ 
                    &\mbox{}& \quad \qquad + \left( \frac{ 1 }{ 2 cp } + \frac{ 1 }{ 2 cq } - \frac{ u + v - tr }{ c^{2} pq } \right) \sum_{ \lambda^{ c } = 1 \not= \lambda } \frac{ \lambda^{ - tr } }{ 1 - \lambda } - \frac{ 1 }{ c^{2} pq } \sum_{ \lambda^{ c } = 1 \not= \lambda } \frac{ \lambda^{ - tr + 1 } }{ ( 1 - \lambda )^{2} } \\ 
                    &\mbox{}& \quad \qquad + \frac{ 1 }{ cp }  \sum_{ \lambda^{ cp } = 1 \not= \lambda^{c}  } \frac{ \lambda^{ v - t r }  }{ \left( 1 - \lambda^{ cq  }  \right) \left( 1 - \lambda \right)  } + \frac{ 1 }{ cq }  \sum_{ \lambda^{ cq } = 1 \not= \lambda^{c}  } \frac{ \lambda^{ u - t r }  }{ \left( 1 - \lambda^{ cp  }  \right) \left( 1 - \lambda \right)  } \ , \end{eqnarray*} 
where $u$ and $v$ are given by (\ref{uandv}). 
\end{proposition} 
It will be useful to have the Laurent expansion of the factors of our generating function. The following lemma will provide a bridge between the residue method and the partial fraction method.  
\begin{lemma}\label{laurent} Let $a,b$ be positive integers, and $ \lambda^{a} = 1 $. Then 
  \[ \frac{ 1 }{ 1 - z^{ ab } } = - \frac{ \lambda }{ ab } ( z - \lambda )^{-1} + \frac{ ab - 1 }{ 2ab } + O ( z - \lambda ) \ . \]
\end{lemma} 
{\it Proof.} First, 
  \[ \mbox{Res} \Bigl( \frac{ 1 }{ 1 - z^{ ab } } , z=\lambda \Bigr) = \lim_{ z \to \lambda } \frac{ z - \lambda }{ 1 - z^{ ab } } = - \frac{ \lambda }{ ab } \ . \] 
For $ ab = 1 $, the statement is trivial, so we may assume $ ab \geq 2 $. Then the constant term of the Laurent series of $ \frac{ 1 }{ 1 - z^{ ab } } $ can be computed as 
  \begin{eqnarray*} &\mbox{}& \lim_{ z \to \lambda } \left( \frac{ 1 }{ 1 - z^{ ab } } + \frac{ \lambda }{ ab ( z - \lambda ) } \right) = \lim_{ z \to \lambda } \frac{ ab ( z - \lambda ) + \lambda \left( 1 - z^{ ab } \right) }{ ab ( z - \lambda ) \left( 1 - z^{ ab } \right) } = \lim_{ z \to \lambda } \frac{ ab - ab \lambda z^{ ab-1 } }{ ab \left( 1 - z^{ ab } - ( z - \lambda )  ab z^{ ab-1 } \right) } \\ 
                    &\mbox{}& \qquad = \lim_{ z \to \lambda } \frac{ - \lambda (ab-1) z^{ ab-2 } }{ - 2ab z^{ ab-1 } - ( z - \lambda ) ab (ab-1) z^{ ab-2 } } = \frac{ ab - 1 }{ 2 ab } \ . \end{eqnarray*} 
\hfill {} $\Box$ 

{\it First proof of Proposition \ref{closure} (partial fractions).} 
To make life easier, we translate the coefficient of $ z^{ t r  } $ of our generating function, which yields the lattice 
point count, to the constant coefficient of the function 
  \begin{equation}\label{partfunction} \frac{ z^{ u+v-tr }  }{ \left( 1 - z^{ cp  } \right) \left( 1 - z^{ cq  } \right) ( 1 - z ) } \ . \end{equation} 
This is a proper rational function because 
  \[ cp \frac{a}{d} + cq \frac{b}{d} \leq r \] 
($ \T \not= \varnothing $!), which implies 
  \[ u+v-tr-cp-cq-1 = \left[ \frac{ ta - 1 }{ d }  \right] cp + \left[ \frac{ tb - 1 }{ d }  \right] cq - tr - 1 < \frac{ta}{d} cp + \frac{tb}{d} cq - tr - 1 \leq -1 \ . \] 
By expanding (\ref{partfunction}) into partial fractions 
  \begin{eqnarray*} &\mbox{}& \frac{ z^{ u+v-tr }  }{ \left( 1 - z^{ cp  } \right) \left( 1 - z^{ cq  } \right) ( 1 - z ) } = \\ 
     &\mbox{}& \qquad \sum_{ \lambda^{cp} = 1 \not= \lambda^{c} }  \frac{ A_{ \lambda } }{ z - \lambda } + \sum_{ \lambda^{cq} = 1 \not= \lambda^{c} } \frac{ B_{ \lambda } }{ z - \lambda } + \sum_{ \lambda^{c} = 1 \not= \lambda } \left( \frac{ C_{ \lambda } }{ z - \lambda } + \frac{ D_{ \lambda } }{ ( z - \lambda )^{2} } \right) + \sum_{k=1}^{3} \frac{ E_{k} }{ (z-1)^{k} } + \sum_{k=1}^{tr-u-v} \frac{ F_{k} }{ z^{k} } \ , \end{eqnarray*} 
we can compute $ L \left( \overline{\T} , t \right) $ as the constant coefficient of the right--hand side: 
  \begin{equation}\label{partial} L \left( \overline{\T} , t \right) = - \sum_{ \lambda^{cp} = 1 \not= \lambda^{c} }  \frac{ A_{ \lambda } }{ \lambda } \ - \sum_{ \lambda^{cq} = 1 \not= \lambda^{c} } \frac{ B_{ \lambda } }{ \lambda } \ + \sum_{ \lambda^{c} = 1 \not= \lambda } \left( - \frac{ C_{ \lambda } }{ \lambda } + \frac{ D_{ \lambda } }{ \lambda^{2} } \right) -  E_{1} + E_{2} - E_{3} \ . \end{equation} 
The computation of the coefficients $ A_{ \lambda } $ for $ \lambda^{cp} = 1 \not= \lambda^{c} $ is straightforward: 
  \[ A_{ \lambda } = \lim_{ z \to \lambda } \frac{ ( z - \lambda ) z^{ u+v-tr }  }{ \left( 1 - z^{ cp  } \right) \left( 1 - z^{ cq  } \right) ( 1 - z ) } = \frac{\lambda^{ v-tr }  }{ \left( 1 - \lambda^{ cq  } \right) ( 1 - \lambda ) } \lim_{ z \to \lambda } \frac{ ( z - \lambda ) }{ 1 - z^{ cp  } } = - \frac{\lambda^{ v-tr+1 }  }{ cp \left( 1 - \lambda^{ cq  } \right) ( 1 - \lambda ) }  \ . \] 
Similarly, we obtain for the $cq$th roots of unity $ \lambda^{cq} = 1 \not= \lambda^{c} $ 
  \[ B_{ \lambda } = - \frac{\lambda^{ u-tr+1 }  }{ cq \left( 1 - \lambda^{ cp  } \right) ( 1 - \lambda ) } \ . \] 
The coefficients $ D_{\lambda} $ and $ C_{\lambda} $ are the two leading coefficients of the Laurent series of 
(\ref{partfunction}) about a nontrivial $c$th root of unity $\lambda$. Using Lemma \ref{laurent}, they are easily seen to be 
  \[ D_{ \lambda } = \frac{\lambda^{ -tr+2 }  }{ c^{2} p q ( 1 - \lambda ) } \] 
and 
  \[ C_{ \lambda } = \left( - \frac{1}{2cp} - \frac{1}{2cq} + \frac{ u+v-tr+1 }{ c^{2} p q } \right) \frac{\lambda^{ -tr+1 }  }{ 1 - \lambda } + \frac{\lambda^{ -tr+2 }  }{ c^{2} p q ( 1 - \lambda )^{2} } \ . \] 
Finally, we obtain the coefficients $ E_{k} $ from the Laurent series of (\ref{partfunction}) about $z=1$ (by hand or, preferably, using a computer algebra system) as 
  \[ E_{3} = - \frac{1}{ c^{2} pq } \ , \qquad E_{2} = - \frac{ u+v-tr+1 }{ c^{2} p q } + \frac{1}{2cp} + \frac{1}{2cq} \ , \]  
and 
  \begin{eqnarray*} &\mbox{}& E_{1} = - \frac{ (u+v-tr)^{2} }{ 2 c^{2} p q } + \frac{ u+v-tr }{ 2 } \left( - \frac{ 1 }{ c^{2} p q } + \frac{1}{cp} + \frac{1}{cq}  \right) + \frac{1}{4} \left( \frac{1}{cp} + \frac{1}{cq} - 1 \right) - \frac{1}{12} \left( \frac{p}{q} + \frac{ 1 }{ c^{2} p q } + \frac{q}{p} \right) \ . \end{eqnarray*} 
Putting these ingredients into (\ref{partial}) yields the statement. 
\hfill {} $\Box$ 

{\it Second proof of Proposition \ref{closure} (residue theorem).} 
The sought--after Taylor coefficient of (\ref{function}) can be shifted to a residue: 
  \begin{equation}\label{residue} L \left( \overline{\T} , t \right) = \mbox{Res} \left( \frac{ z^{ u + v - t r  - 1 } }{ \left( 1 - z^{ cp }  \right) \left( 1 - z^{ cq }  \right) \left( 1 - z \right) } , z=0 \right) \end{equation}
If the right--hand side of (\ref{residue}) counts the lattice points in $ t \T $, then what we 
have to do is compute the other residues of 
  \[ f (z) := \frac{ z^{ u + v - t r  - 1 } }{ \left( 1 - z^{ cp }  \right) \left( 1 - z^{ cq }  \right) \left( 1 - z \right) } \ , \]
and use the residue theorem for the sphere $\C \cup \left\{ \infty \right\} $. Aside from 0, $ f $ has poles at 
all $ cp $th and $ cq $th roots of unity; note that the nonemptyness of $ \T $ implies Res($ f(z), z = \infty $) = 0.  

The residue at $ z=1 $ can be easily calculated as 
  \begin{eqnarray*} &\mbox{}& \mbox{Res} \Bigl( f (z) , z=1 \Bigr) = \mbox{Res} \Bigl( e^{ z }  f (e^{ z } ) , z=0 \Bigr) = - \frac{ 1 }{ 2 c^{2} p q  } \left( u + v - t r \right)^{ 2 }  \\ 
                    &\mbox{}& \quad + \frac{ 1 }{ 2 } \left( u + v - t r \right) \left( \frac{ 1 }{ cp  } + \frac{ 1 }{ cq  } + \frac{ 1 }{ c^{2} p q  }  \right) - \frac{ 1 }{ 4 } \left( 1 + \frac{ 1 }{ cp  } + \frac{ 1 }{ cq  }  \right) - \frac{ 1 }{ 12 } \left( \frac{ p  }{ q  } + \frac{ q  }{ p  } + \frac{ 1 }{ c^{2} p q  }  \right) \ . \end{eqnarray*} 
It remains to compute the residues at the nontrivial roots of unity. Let $ \lambda $ be a nontrivial $c$th roots of unity. 
Putting the Laurent expansions of the different factors of $f$ together, the residue of $f$ at $ \lambda $ can be easily derived via Lemma \ref{laurent}  as 
  \[ \mbox{Res} \Bigl( f(z) , z=\lambda \Bigr) = \left( - \frac{ u + v - tr }{ c^{2} pq } + \frac{ 1 }{ 2 cp } + \frac{ 1 }{ 2 cq } \right) \frac{ \lambda^{ - tr } }{ 1 - \lambda } + \frac{ \lambda^{ - tr + 1 } }{ c^{2} pq ( 1 - \lambda )^{2} } \ . \] 
Note that we have to add this expression over all nontrivial $c$th roots of unity. 

Now let $ \lambda^{ cp } = 1 \not= \lambda^{c} $. Since $p$ and $q$ are relatively prime, $f$ has a simple pole at $ \lambda $, whose residue can be determined easily using Lemma \ref{laurent}: 
  \[ \mbox{Res} \Bigl( f (z) , z=\lambda \Bigr) = \frac{ \lambda^{ v - t r - 1}  }{ \left( 1 - \lambda^{ cq  }  \right) \left( 1 - \lambda \right)  } \ \mbox{Res} \left( \frac{ 1 }{ 1 - \lambda^{ cp  }  }  , z=\lambda \right) = - \frac{ \lambda^{ v - t r }  }{ cp \left( 1 - \lambda^{ cq  }  \right) \left( 1 - \lambda \right)  } \ . \] 
Again, we have to add up all these $\lambda$'s. Finally, we obtain a similar expression for the $cq$th roots of unity, and 
the statement of the proposition follows by rewriting (\ref{residue}) by means of the residue theorem. 
\hfill {} $\Box$ 

In the following section, we will further describe the finite sums appearing in the lattice point count operators; consequently, we will be able to make statements about their computational complexity. 


\section{Using the Dedekind--Rademacher sums as building blocks}\label{sums} 
We will now take a closer look at the finite sums over roots of unity appearing in Proposition 
\ref{closure}, namely, 
  \[ \frac{ 1 }{ cp }  \sum_{ \lambda^{ cp } = 1 \not= \lambda^{c}  } \frac{ \lambda^{ w }  }{ \left( 1 - \lambda^{ cq  }  \right) \left( 1 - \lambda \right)  }  \]
for some integers $c,p,q,w$, where $p$ and $q$ are relatively prime.  
The fact that this is a finite Fourier series in $w$ and the appearance of two factors in the 
denominator suggest the use of the well--known Convolution Theorem for finite Fourier series: 
\begin{theorem}\label{convolution} Let $ \displaystyle f(t) = \frac{ 1 }{ N }  \sum_{ \lambda^{N} = 1 } a_{ \lambda } \lambda^{ t } $ and 
$ \displaystyle g(t) = \frac{ 1 }{ N } \sum_{ \lambda^{N} = 1 } b_{ \lambda } \lambda^{ t } $. Then 
  \[ \frac{ 1 }{ N } \sum_{ \lambda^{N} = 1 } a_{ \lambda } b_{ \lambda } \lambda^{ t } = \sum_{ m=0 }^{ N-1 } f(t-m) g(m) \ .  \] 
\hfill {} $\Box$ 
\end{theorem} 
The key ingredient to be able to apply this theorem to our case is  
\begin{lemma}\label{key} For $ p, t \in \Z $, 
  \[ \frac{1}{p} \sum_{ \lambda^{ p } = 1 \not= \lambda } \frac{ \lambda^{ t } }{ \lambda - 1 } = \left( \left( \frac{-t}{p} \right) \right) + \frac{1}{2p} \ . \]
\end{lemma} 
Recall that $ (( x )) = x - [ x ] - 1/2 $. This lemma is well--known (see, e.g., \cite{grosswald}, p.~14), 
however, for sake of completeness we give a short proof based on the methods of section \ref{mattSec}: 

{\it Proof.} Consider the interval $ \I := [ 0 , \frac{ 1 }{ p } ] $, viewed as a one--dimensional 
polytope. Then the lattice point count in the dilated interval is clearly 
  \begin{equation}\label{inttrivial} L \left( \overline{\I} , t \right) = \left[ \frac{ t }{ p }  \right] + 1 \ . \end{equation} 
On the other hand, we can write this number, by applying the ideas in section \ref{mattSec}, as 
  \[ L \left( \overline{\I} , t \right) = \mbox{Res} \left( \frac{ z^{ - t - 1 } }{ \left( 1 - z^{ p }  \right) \left( 1 - z \right) } , z=0 \right) \ . \] 
(Equivalently, we could expand this generating function into partial fractions.) 
Using again the residue theorem, this can be rewritten as 
  \begin{equation}\label{intresidue} L \left( \overline{\I} , t \right) = \frac{ t }{ p } + \frac{ 1 }{ 2p } + \frac{ 1 }{ 2 } - \frac{1}{p} \sum_{ \lambda^{ p } = 1 \not= \lambda } \frac{ \lambda^{ -t } }{  \lambda - 1 } \ . \end{equation} 
Comparing (\ref{inttrivial}) with (\ref{intresidue}) yields the statement. 
\hfill {} $\Box$ 
\begin{corollary}\label{keycor} For $ c, p, q, t \in \Z , (p,q) = 1 $, 
  \[ \frac{1}{cp} \sum_{ \lambda^{ cp } = 1 \not= \lambda^{c} } \frac{ \lambda^{ t } }{ 1 - \lambda^{cq} } = 
       \left\{ \begin{array}{ll} - \left( \left( \frac{- q^{-1} t }{cp} \right) \right) - \frac{1}{2p} & \mbox{ if } c | t \\ 
                                 0 & \mbox{ otherwise. } \end{array} \right. \]
Here, $ q q^{-1} \equiv 1 $ mod $p$. 
\end{corollary} 
{\it Proof.} If $ c|t $, write $ t=cw $ to obtain 
  \begin{eqnarray*} &\mbox{}& \frac{1}{cp} \sum_{ \lambda^{ cp } = 1 \not= \lambda^{c} } \frac{ \lambda^{ t } }{ 1 - \lambda^{cq} } = \frac{1}{cp} \sum_{ \lambda^{ cp } = 1 \not= \lambda^{c} } \frac{ \lambda^{ cw } }{ 1 - \lambda^{cq} } = \frac{1}{p} \sum_{ \lambda^{ p } = 1 \not= \lambda } \frac{ \lambda^{ w } }{ 1 - \lambda^{q} } = \frac{1}{p} \sum_{ \lambda^{ p } = 1 \not= \lambda } \frac{ \lambda^{ q^{-1} w } }{ 1 - \lambda } \\ 
                    &\mbox{}& \qquad \stackrel{ (\star) }{ = } - \left( \left( \frac{- q^{-1} w }{p} \right) \right) - \frac{1}{2p} = - \left( \left( \frac{- q^{-1} t }{cp} \right) \right) - \frac{1}{2p} \ . \end{eqnarray*} 
Here, $(\star)$ follows from Lemma \ref{key}. If $c$ does not divide $t$, let $ \xi = e^{ 2 \pi i / cp } $. Then 
  \begin{eqnarray*} &\mbox{}& \frac{1}{cp} \sum_{ \lambda^{ cp } = 1 \not= \lambda^{c} } \frac{ \lambda^{ t } }{ 1 - \lambda^{cq} } = \frac{1}{cp} \sum_{ m=1 }^{ p-1 } \sum_{ n=0 }^{ c-1 } \frac{ \xi^{ (mc+np) t } }{ 1 - \xi^{ (mc+np) cq } } = \frac{1}{cp} \sum_{ n=0 }^{ c-1 } \xi^{ npt } \sum_{ m=1 }^{ p-1 } \frac{ \xi^{ mct } }{ 1 - \xi^{ m c^{2} q } } = 0 \ . \end{eqnarray*} 
\hfill {} $\Box$ 
\begin{corollary}\label{raddedsum}  For $ c, p, q, t \in \Z , (p,q) = 1 $, 
  \[ \frac{1}{cp} \sum_{ \lambda^{ cp } = 1 \not= \lambda^{c} } \frac{ \lambda^{ -t } }{ \left( 1 - \lambda^{cq} \right)  \left( 1 - \lambda \right) } = - \sigma \left( q,p, \frac{ t }{ c } \right) - \left( \left( \frac{ t }{ cp } \right) \right) + \frac{1}{2p} \left( \left( \frac{ t }{ c } \right) \right) \ . \] \end{corollary} 
{\it Proof.} We will repeatedly use the periodicity of the sawtooth function. One consequence is, for $ p \in \Z, x \in \R $,  
  \begin{equation}\label{trivial} \sum_{ m=0 }^{ p-1 } \left( \left( \frac{ m+x }{ p } \right) \right) = ((x)) \ , \end{equation} 
the proof of which is left as an exercise (\cite{grosswald}, p.~4). 
Now by Lemma \ref{key}, 
  \begin{eqnarray*} &\mbox{}& \frac{1}{cp} \sum_{ \lambda^{ cp } = 1 \not= \lambda^{c} } \frac{ \lambda^{ t } }{ \left( 1 - \lambda \right) } = \frac{1}{cp} \sum_{ \lambda^{ cp } = 1 \not= \lambda } \frac{ \lambda^{ t } }{ \left( 1 - \lambda \right) } - \frac{1}{cp} \sum_{ \lambda^{ c } = 1 \not= \lambda } \frac{ \lambda^{ t } }{ \left( 1 - \lambda \right) } \\ 
                   &\mbox{}& \qquad = - \left( \left( \frac{-t}{cp} \right) \right) - \frac{1}{2cp} - \frac{1}{p} \left( - \left( \left( \frac{-t}{c} \right) \right) - \frac{1}{2c}  \right) = - \left( \left( \frac{-t}{cp} \right) \right) + \frac{1}{p} \left( \left( \frac{-t}{c} \right) \right) \ . \end{eqnarray*} 
Finally we use the Convolution Theorem \ref{convolution} and Corollary \ref{keycor} to obtain 
  \begin{eqnarray*} &\mbox{}& \frac{1}{cp} \sum_{ \lambda^{ cp } = 1 \not= \lambda^{c} } \frac{ \lambda^{ t } }{ \left( 1 - \lambda^{cq} \right)  \left( 1 - \lambda \right) } = \\ 
                   &\mbox{}& \qquad = \sum_{ {m=0} \atop {c|m} }^{ cp-1 } \left( - \left( \left( \frac{ - q^{-1} m }{ cp } \right) \right) - \frac{1}{2p} \right) \left( - \left( \left( \frac{ - (t-m) }{ cp } \right) \right) + \frac{1}{p} \left( \left( \frac{ - (t-m) }{ c } \right) \right)  \right) \\ 
                   &\mbox{}& \qquad = \sum_{ {k=0} }^{ p-1 } \left( \left( \frac{ - q^{-1} k }{ p } \right) \right)  \left( \left( \frac{ k }{ p } - \frac{ t }{ cp } \right) \right) - \frac{1}{p} \sum_{ {k=0} }^{ p-1 } \left( \left( \frac{ - q^{-1} k }{ p } \right) \right)  \left( \left( \frac{ -t }{ c } \right) \right) \\ 
                   &\mbox{}& \qquad \quad + \frac{1}{2p} \sum_{ {k=0} }^{ p-1 } \left( \left( \frac{ k }{ p } - \frac{ t }{ cp } \right) \right) -  \frac{1}{2p^{2}} \sum_{ {k=0} }^{ p-1 } \left( \left( \frac{ -t }{ c } \right) \right) \\
                   &\mbox{}& \qquad \stackrel{ (\ref{trivial}) }{ = } \sum_{ {k=0} }^{ p-1 } \left( \left( \frac{ - k }{ p } \right) \right)  \left( \left( \frac{ qk }{ p } - \frac{ t }{ cp } \right) \right) + \frac{1}{2p} \left( \left( \frac{ -t }{ c } \right) \right) + \frac{1}{2p} \left( \left( \frac{ -t }{ c } \right) \right) - \frac{1}{2p} \left( \left( \frac{ -t }{ c } \right) \right) \\ 
                   &\mbox{}& \qquad \stackrel{ (\ref{trivial}) }{ = } - \sum_{ {k=0} }^{ p-1 } \left( \left( \frac{ k }{ p } \right) \right)  \left( \left( \frac{ qk }{ p } - \frac{ t }{ cp } \right) \right) - \left( \left( \frac{ -t }{ cp } \right) \right) + \frac{1}{2p} \left( \left( \frac{ -t }{ c } \right) \right) \ . \end{eqnarray*}
The statement follows now by definition of the Dedekind--Rademacher sum (\ref{s_t}). 
\hfill {} $\Box$ 

This corollary describes the finite sums in Proposition \ref{closure}. One of them actually turns out to be of an even simpler form. To show this, we first need to rewrite Proposition \ref{closure} for the special case where $\T$ has the origin as a vertex: 
\begin{proposition}\label{special} For the rectangular rational triangle $\T$ given by (\ref{triangle}) with $ a=b=0 $, $ c=r=1 $, and $p$ and $q$ relatively prime,  
  \begin{eqnarray*} &\mbox{}& L \left( \overline{\T} , t \right) = \frac{ t^{ 2 } }{ 2 p q  } + \frac{ t }{ 2 } \left( \frac{ 1 }{ p  } + \frac{ 1 }{ q  } + \frac{ 1 }{ p q  }  \right) + \frac{ 1 }{ 4 } + \frac{ 1 }{ 12 } \left( \frac{ p  }{ q  } + \frac{ q  }{ p  } + \frac{ 1 }{ p q  }  \right) \\ 
                    &\mbox{}& \quad \qquad - \sigma \left( q, p, t \right) - \sigma \left( p, q, t \right) - \left( \left( \frac{ t }{ p } \right) \right) - \left( \left( \frac{ t }{ q } \right) \right)  \ . \end{eqnarray*} \end{proposition} 
{\it Proof.} Proposition \ref{closure} gives for this special case
  \begin{eqnarray*} &\mbox{}& L \left( \overline{\T} , t \right) = \frac{ t^{2} }{ 2 p q  } + \frac{ t }{ 2 } \left( \frac{ 1 }{ p  } + \frac{ 1 }{ q  } + \frac{ 1 }{ p q  }  \right) + \frac{ 1 }{ 4 } \left( 1 + \frac{ 1 }{ p  } + \frac{ 1 }{ q  }  \right) + \frac{ 1 }{ 12 } \left( \frac{ p  }{ q  } + \frac{ q  }{ p  } + \frac{ 1 }{ p q  }  \right) \\ 
                    &\mbox{}& \quad \qquad + \frac{ 1 }{ p }  \sum_{ \lambda^{ p } = 1 \not= \lambda } \frac{ \lambda^{ - t }  }{ \left( 1 - \lambda^{ q  }  \right) \left( 1 - \lambda \right)  } + \frac{ 1 }{ q }  \sum_{ \mu^{ q } = 1 \not= \mu } \frac{ \mu^{ - t }  }{ \left( 1 - \mu^{ p  }  \right) \left( 1 - \mu \right)  } \ . \end{eqnarray*} 
The statement follows now with Corollary \ref{raddedsum}. 
\hfill {} $\Box$ 

We use this Proposition to show
\begin{lemma}\label{trivdedsum} For $ p, t \in \Z $, 
\[ \sigma \left( 1, p, t \right) = \sum_{ k=0 }^{ p-1 } \left( \left( \frac{ k+t }{ p } \right) \right) \left( \left( \frac{ k }{ p } \right) \right) = - \frac{p}{24} + \frac{1}{6p} + \frac{p}{2} \left( \left( \frac{ t }{ p } \right) \right)^{2}   \ . \] \end{lemma} 
{\it Proof.} Consider the triangle $ \Delta := \left\{ (x,y) \in \R^{2} \ : \ x + py \leq 1 \right\} $ and its integer dilates. By summing over vertical line segments in the triangle, we obtain  
  \begin{eqnarray} &\mbox{}& L \left( \Delta , t \right) = \sum_{ m=0 }^{ \left[ \frac{ t }{ p } \right] } ( t - pm + 1 ) = (t+1) \left( \left[ \frac{ t }{ p } \right] + 1 \right) - \frac{p}{2} \left[ \frac{ t }{ p } \right] \left( \left[ \frac{ t }{ p } \right] + 1 \right) \nonumber \\ 
                    &\mbox{}& \qquad = \frac{ t^{2} }{ 2p } + \left( \frac{1}{p} + \frac{1}{2} \right) t + \frac{1}{2} + \frac{p}{8} - \left( \left( \frac{ t }{ p } \right) \right) - \frac{p}{2} \left( \left( \frac{ t }{ p } \right) \right)^{2} \ . \label{once} \end{eqnarray} 
On the other hand, we can compute the same number via Proposition \ref{special}: 
  \begin{equation}\label{twice} L \left( \Delta , t \right) =  \frac{ t^{2} }{ 2p } + \frac{t}{2} \left( \frac{2}{p} + 1 \right) + \frac{1}{4} + \frac{ 1 }{ 12 } \left( p + \frac{ 2 }{ p  } \right) - \sigma \left( 1, p, t \right) - \frac{1}{4} - \left( \left( \frac{ t }{ p } \right) \right) + \frac{1}{2} \ . \end{equation} 
Again we used (\ref{trivial}). Equating (\ref{once}) with (\ref{twice}) yields the statement. 
\hfill {} $\Box$ 

Using these ingredients, we can finally restate Proposition \ref{closure} as our main theorem:  
\begin{theorem}\label{bigthm} For the rectangular rational triangle $\T$ given by (\ref{triangle}), 
  \begin{eqnarray*} &\mbox{}& L \left( \overline{\T} , t \right) = \frac{ 1 }{ 2 c^{2} p q  } ( tr-u-v )^{2} + ( tr-u-v ) \left( \frac{ 1 }{ 2cp } + \frac{ 1 }{ 2cq } + \frac{ 1 }{ c^{2} p q  } + \frac{ 1 }{ cpq } \left( \left( \frac{ tr }{ c } \right) \right) \right) \\ 
                    &\mbox{}& \quad \qquad + \frac{ 1 }{ 4 } + \frac{ 1 }{ 12 } \left( \frac{ p }{ q } + \frac{ q }{ p } \right) - \frac{ 1 }{ 24 pq } + \frac{ 1 }{ c^{2} p q } + \left( \frac{1}{2p} + \frac{1}{2q} \right) \left( \left( \frac{ tr }{ c } \right) \right)  \\
                    &\mbox{}& \quad \qquad - \left( \left( \frac{ tr-v }{ cp } \right) \right) - \left( \left( \frac{ tr-u }{ cq } \right) \right) + \left( \frac{ 1 }{ cpq } - \frac{ 1 }{ 2p } - \frac{ 1 }{ 2q } \right) \left( \left( \frac{ tr }{ c } \right) \right) \\ 
                    &\mbox{}& \quad \qquad + \frac{ 1 }{ cpq } \left( \left( \frac{ tr-1 }{ c } \right) \right) + \frac{ 1 }{ 2 pq } \left( \left( \frac{ tr-1 }{ c } \right) \right)^{2} - \sigma \left( q, p, \frac{ tr-v }{ c } \right) - \sigma \left( p, q, \frac{ tr-u }{ c } \right)  \ . \end{eqnarray*} 
Here $u$ and $v$ are given by (\ref{uandv}). 
\end{theorem} 
{\it Proof.} By Lemma \ref{key}, 
  \begin{equation}\label{simp1} \frac{1}{c} \sum_{ \lambda^{ c } = 1 \not= \lambda } \frac{ \lambda^{ w } }{ 1 - \lambda } = - \left( \left( \frac{-w}{c} \right) \right) - \frac{1}{2c} \ . \end{equation} 
By Corollary \ref{raddedsum} and Lemma \ref{trivdedsum}, 
  \begin{eqnarray} &\mbox{}& \frac{1}{c} \sum_{ \lambda^{ c } = 1 \not= \lambda } \frac{ \lambda^{ w } }{ ( 1 - \lambda )^{2} } = - \sigma \left( 1, c -w \right) - \left( \left( \frac{-w}{c} \right) \right) - \frac{1}{4c} \nonumber \\ 
                   &\mbox{}& \qquad = \frac{c}{24} - \frac{5}{12c} - \left( \left( \frac{ -w }{ c } \right) \right) - \frac{c}{2} \left( \left( \frac{ -w }{ c } \right) \right)^{2} \ . \label{simp2} \end{eqnarray} 
Now simplify the identity in Proposition \ref{closure} by means of (\ref{simp1}), (\ref{simp2}), and Corollary \ref{raddedsum}. 
\hfill {} $\Box$ 

We mention the work of Brion and Vergne \cite{brion} on rational polytopes, 
where they give formulas for the lattice point enumerator of a convex 
rational polytope in terms of certain Todd differential operators. 
These are interesting connections to topology. 
The salient difference in approaches lies in our 
explicit use of Dedekind--Rademacher sums and their computational efficiency. 
Although \cite{brion} offer an approach which might 
eventually yield similar results, the analysis of the equivalence 
between our different approaches would be the content of another 
paper.  The results in \cite[main result on p.~831]{brion} are 
coordinate--free, and have their theoretical merits.  In contrast our results 
are coordinatized and are useful in the context of computer--related applications. 


\section{Remarks and consequences}\label{consequences}  
Rademacher's original definition \cite{rademacher} of his generalization of the Dedekind sum is 
  \[ S (a,b;x,y) := \frac{ 1 }{ b } \sum_{ k=0 }^{ b-1 } \left( \left( \frac{ a (k+y) }{ b } + x \right)  \right)^{ \star }  \left( \left( \frac{ k+y }{ b }  \right)  \right)^{ \star } \ , \]
defined for $ a, b \in \Z , x, y \in \R $. Here, 
  \[ ((x))^{ \star } := \left\{ \begin{array}{ll} ((x)) & \mbox{ if } x \not\in \Z \\ 
                                                  0     & \mbox{ if } x \in \Z \end{array} \right. \] 
is the sawtooth function which also appears in the classical Dedekind sum $ S (a,b;0,0) $. 
The impact of the slightly different definition of the sawtooth function is not crucial for our lattice point count formula. In fact, it is easy to see that 
  \begin{equation}\label{conv} \sigma (a,b,t) = S \left( a,b; \frac t b ,0 \right) - \frac{1}{2} \left( \left( \frac t b \right) \right) \ . \end{equation} 
An important property of $ S (a,b;x,y) $ is Rademacher's reciprocity law \cite{rademacher} 
  \begin{eqnarray*} &\mbox{}& S (a,b;x,y) + S (b,a;y,x) = \\ 
                    &\mbox{}& \qquad = \left\{ \begin{array}{ll} - \frac{ 1 }{ 4 } + \frac{ 1 }{ 12 } \left( \frac{ a }{ b } + \frac{ 1 }{ ab } + \frac{ b }{ a }  \right) & \mbox{ if both } x,y \in \Z \\ 
                                                                 ((x))^{ \star } ((y))^{ \star } + \frac{ 1 }{ 2 } \left( \frac{ a }{ b } \psi_{ 2 } (y) + \frac{ 1 }{ ab } \psi_{ 2 } (ay+bx) + \frac{ b }{ a } \psi_{ 2 } (x) \right) & \mbox{ otherwise, } \end{array} \right. \end{eqnarray*} 
where 
  \[ \psi_{ 2 } (x) := ( x - [x] )^{ 2 } - ( x - [x] ) + \frac{ 1 }{ 6 }   \]
is the periodic second Bernoulli polynomial. 
The equivalent reciprocity law for $ \sigma (a,b,t) $ was first presented in \cite{knuth}; 
we will rediscover it in the last section.  
 
Among other things, Rademacher's reciprocity law allows us to compute $ S (a,b;x,y) $ (and hence $ \sigma (a,b,t) $, the nontrivial 
part of our lattice point count formulas) in polynomial time, by means of a Euclidean algorithm using the first two variables: 
simply note that we can replace $a$ in $ S (a,b;x,y) $ by the least residue of $a$ modulo $b$. 
It is amusing to note that $ \sigma ( a , b, t ) $ appears in the multiplier system of a weight--0 modular form \cite{robins}. 

To complete the picture for an {\it arbitrary} two--dimensional rational polytope $\P$, we return 
to the statements in the introduction. After triangulating $\P$, the problem reduces to rational 
rectangles and the rectangular triangles which were treated above. A lattice point count formula 
for a rational rectangle $\RR$ is easy to obtain: suppose $\RR$ has vertices 
$ \left( \frac{ a_{ 1 }  }{ d } , \frac{ a_{ 2 }  }{ d }  \right) , 
  \left( \frac{ b_{ 1 }  }{ d } , \frac{ a_{ 2 }  }{ d }  \right) , 
  \left( \frac{ b_{ 1 }  }{ d } , \frac{ b_{ 2 }  }{ d }  \right) , 
  \left( \frac{ a_{ 1 }  }{ d } , \frac{ b_{ 2 }  }{ d }  \right) $, 
with $ a_{ j } < b_{ j } $, then it is not hard to see that 
  \[ L \left( \overline{\RR} , t \right) = \left( \left[ \frac{ t b_{ 1 }  }{ d }  \right] - \left[ \frac{ t a_{ 1 } - 1 }{ d }  \right] \right) \left( \left[ \frac{ t b_{ 2 }  }{ d }  \right] - \left[ \frac{ t a_{ 2 } - 1 }{ d }  \right] \right) \ . \]
Together with the above established remark on computability of $ \sigma (a,b,t) $, we summarize 
some statements in 
\begin{theorem}\label{comp} Let $\P$ be a two--dimensional rational polytope. The coefficients of $ L \left( \overline{\P} , t \right) $ 
can be written in terms of the sawtooth function $((\dots))$ and the 
Dedekind--Radema\-cher sum $ \sigma (a,b,t) $. Consequently, the formula given by Theorem \ref{bigthm} for the lattice point count operator can 
be computed in polynomial time in the logarithm of the denominators of the vertices of $ \P $.  
\hfill {} $\Box$ 
\end{theorem} 
Barvinok \cite{barvinok} showed that for any fixed dimension the lattice point enumerator 
of a rational polytope can be computed in polynomial time. The distinction here is that 
we get a simple {\it formula}, which happens to be polynomial--time computable.  

In dimensions greater than 2, things get more complicated. We can still get formulas through 
the methods introduced here; however, even the existence of (possible) three--term reciprocity laws 
for functions appearing in the lattice point count does not guarantee polynomial--time computability. 
The details will be described in a forthcoming paper \cite{bdr}. 


\section{Reciprocity laws}\label{laws} 
As another remark, we can recover the reciprocity law for the classical De\-de\-kind sum (\cite{dedekind}, Chapter 2 of \cite{grosswald}) from our formulas: 
\begin{corollary}[Dedekind]\label{classic} 
  \[ S(a,b;0,0) + S(b,a;0,0) = - \frac{1}{4} + \frac{1}{12} \left( \frac{a}{b} + \frac{1}{ab} + \frac{b}{a} \right) \ . \] 
\end{corollary} 
{\it Proof.} It is well known \cite{ehrhart} that the constant term in the Ehrhart polynomial, the integer--point enumeration function 
of a {\it lattice polytope} (that is, a polytope with integral vertices) equals the Euler characteristic of the polytope. 
Consider the simplest case of our triangle mentioned in Proposition \ref{special}. If we dilate this polytope by $t=pqw$, that is, only by multiples of $pq$, 
we obtain the dilates of a lattice polytope $\P$. Proposition \ref{special} simplifies for these $t$ to 
  \begin{eqnarray*} &\mbox{}& L \left( \overline{\P} , w \right) = \frac{ pq w^{ 2 } }{ 2 } + \frac{ w }{ 2 } \left( p + q + 1 \right) + \frac{ 1 }{ 4 } + \frac{ 1 }{ 12 } \left( \frac{ p  }{ q  } + \frac{ q  }{ p  } + \frac{ 1 }{ p q  }  \right) - \sigma \bigl( q, p, 0 \bigr) - \sigma \bigl( p, q, 0 \bigr) + 1 \ . \end{eqnarray*} 
On the other hand, we know that the constant term is the Euler characteristic of $\P$ and hence equals 1, which yields the identity 
  \[ \frac{ 1 }{ 4 } + \frac{ 1 }{ 12 } \left( \frac{ p  }{ q  } + \frac{ q  }{ p  } + \frac{ 1 }{ p q  }  \right) - \sigma \bigl( q, p, 0 \bigr) - \sigma \bigl( p, q, 0 \bigr) = 0 \ . \] 
The statement follows now by rewriting the Dedekind sums in terms of the original definition via (\ref{conv}). 
\hfill {} $\Box$ 

As a concluding consequence of our formulas, we rederive a reciprocity law due to Gessel \cite{gessel}, 
interpreting it at the same time geometrically. 
\begin{corollary}[Gessel]\label{ges} Let $p$ and $q$ be relatively prime and suppose that $t$ is an integer such that $ 1 \leq t \leq p+q $. Then 
  \begin{eqnarray*} &\mbox{}& \frac{ 1 }{ p }  \sum_{ \lambda^{ p } = 1 \not= \lambda } \frac{ \lambda^{ t }  }{ \left( 1 - \lambda^{ q } \right) \left( 1 - \lambda \right)  } +  \frac{ 1 }{ q }  \sum_{ \lambda^{ q } = 1 \not= \lambda } \frac{ \lambda^{ t }  }{ \left( 1 - \lambda^{ p  }  \right) \left( 1 - \lambda \right)  } \\ 
                    &\mbox{}& \qquad = - \frac{ t^{2} }{ 2pq } + \frac{t}{2} \left( \frac{1}{p} + \frac{1}{q} + \frac{1}{pq} \right) - \frac{1}{4} \left( \frac{1}{p} + \frac{1}{q} + 1 \right)  - \frac{1}{12} \left( \frac{p}{q} + \frac{1}{pq} + \frac{q}{p} \right) \ . \end{eqnarray*} 
\end{corollary} 
It is easy to see that the reciprocity law for classical Dedekind sums (Corollary \ref{classic}) is a special 
case of Gessel's theorem. Before proving Gessel's theorem below,  
we find it useful to have the lattice point count operator for the interior of our polytope. 
The following central theorem, conjectured by Ehrhart \cite{ehrhart} and first proved by 
Macdonald \cite{macdonald}, enables us to do this: 
\begin{theorem}[Ehrhart--Macdonald Reciprocity Law]\label{reclaw} If a rational polytope $\P$  
is homeomorphic to a $d$--sphere then  
  \[ L ( {\cal P}^{\circ}  , -t ) = (-1)^{ d } L ( \overline{\cal P} , t ) \ . \]
\hfill {} $\Box$ 
\end{theorem} 
Note that this theorem allows us to conclude a computability statement for the lattice point count in the 
{\it interior} of a two--dimensional rational polytope similar to Theorem \ref{comp}. 
Using Theorem \ref{reclaw} we get from Proposition \ref{closure} the 
\begin{corollary}\label{interior} For the rectangular rational triangle $\T$ given by (\ref{triangle}) with $ a=b=0 $, $ c=r=1 $, and $p$ and $q$ relatively prime,  
  \begin{eqnarray*} &\mbox{}& L \left( \T^{\circ} , t \right) = \frac{ t^{2} }{ 2 p q  } - \frac{ t }{ 2 } \left( \frac{ 1 }{ p } + \frac{ 1 }{ q } + \frac{ 1 }{ p q  }  \right) + \frac{ 1 }{ 4 } \left( 1 + \frac{ 1 }{ p  } + \frac{ 1 }{ q  }  \right) + \frac{ 1 }{ 12 } \left( \frac{ p  }{ q  } + \frac{ q  }{ p  } + \frac{ 1 }{ p q  }  \right) \\ 
                    &\mbox{}& \quad \qquad + \frac{ 1 }{ p }  \sum_{ \lambda^{ p } = 1 \not= \lambda } \frac{ \lambda^{ t }  }{ \left( 1 - \lambda^{ q  }  \right) \left( 1 - \lambda \right)  } + \frac{ 1 }{ q }  \sum_{ \lambda^{ q } = 1 \not= \lambda } \frac{ \lambda^{ t }  }{ \left( 1 - \lambda^{ p  }  \right) \left( 1 - \lambda \right)  } \ . \end{eqnarray*} 
\hfill {} $\Box$ 
\end{corollary} 
Alternatively, we could have derived Corollary \ref{interior} from scratch in a similar way as Proposition \ref{closure}. 

{\it Proof of Corollary \ref{ges}.} Consider dilates of the triangle given in Corollary \ref{interior}, that is, 
  \[ t \T^{\circ} = \left\{ (x,y) \in \R^{ 2 } : \ x, y > 0 , \ p x + q y < t \right\} \ . \] 
By the very definition, $ t \T^{\circ} $ does not contain any integer points for $ 1 \leq t \leq p+q $, in other words, 
$ L \left( \T^{\circ} , t \right) = 0 $. Hence Corollary \ref{interior} yields an identity for these values of $t$: 
  \begin{eqnarray*} &\mbox{}& 0 = \frac{ t^{2} }{ 2 p q  } - \frac{ t }{ 2 } \left( \frac{ 1 }{ p } + \frac{ 1 }{ q } + \frac{ 1 }{ p q  }  \right) + \frac{ 1 }{ 4 } \left( 1 + \frac{ 1 }{ p  } + \frac{ 1 }{ q  }  \right) + \frac{ 1 }{ 12 } \left( \frac{ p  }{ q  } + \frac{ q  }{ p  } + \frac{ 1 }{ p q  }  \right) \\ 
                    &\mbox{}& \quad \qquad + \frac{ 1 }{ p }  \sum_{ \lambda^{ p } = 1 \not= \lambda } \frac{ \lambda^{ t }  }{ \left( 1 - \lambda^{ q  }  \right) \left( 1 - \lambda \right)  } + \frac{ 1 }{ q }  \sum_{ \lambda^{ q } = 1 \not= \lambda } \frac{ \lambda^{ t }  }{ \left( 1 - \lambda^{ p  }  \right) \left( 1 - \lambda \right)  } \ . \end{eqnarray*} 
 
We can rephrase Corollary \ref{ges} in terms of Dedekind--Rademacher sums by means of Corollary \ref{raddedsum}:  
\begin{corollary}\label{gesvar} Let $p$ and $q$ be relatively prime and suppose that $t$ is an integer such that $ 1 \leq t \leq p+q $. Then  
  \begin{eqnarray*} &\mbox{}& \sigma \bigl( q, p, -t \bigr) + \sigma \bigl( p, q, -t \bigr) \stackrel{ \mbox{\rm \footnotesize def} }{ = } \ \sum_{ k=0 }^{ p-1 } \left( \left( \frac{ q k - t }{ p } \right) \right) \left( \left( \frac{ k }{ p } \right) \right) + \sum_{ k=0 }^{ q-1 } \left( \left( \frac{ p k - t }{ q } \right) \right) \left( \left( \frac{ k }{ q } \right) \right) \\  
                    &\mbox{}& \qquad = \frac{ t^{2} }{ 2pq } - \frac{t}{2} \left( \frac{1}{p} + \frac{1}{q} + \frac{1}{pq} \right) + \frac{1}{4} + \frac{1}{12} \left( \frac{p}{q} + \frac{1}{pq} + \frac{q}{p} \right) - \left( \left( \frac{ -t }{ p } \right) \right) - \left( \left( \frac{ -t }{ q } \right) \right)  \ . \end{eqnarray*}  
\hfill {} $\Box$  
\end{corollary}  
This version of Gessel's reciprocity law bears something surprising: its form is essentially identical to Knuth's  
version \cite{knuth} of Rademacher's reciprocity law for the Dedekind--Rademacher sums, with two differences:  
Gessel's theorem requires $t$ to be an integer, whereas $ t \in \R $ in Knuth's reciprocity law.  
On the other hand, the conditions in Knuth's theorem are $ p < q $ and $ 0 \leq t \leq q $,  
which suffices for all practical purposes, however, Gessel's range on $t$ is bigger.  
In fact, is not hard to remove the integrality condition on $t$, which unifies Knuth's and Gessel's reciprocity theorem:  
\begin{theorem} Let $p$ and $q$ be relatively prime and suppose that $t$ is a real number such that $ 1 \leq t \leq p+q $. Then  
  \begin{eqnarray*} &\mbox{}& \sigma \bigl( q, p, -t \bigr) + \sigma \bigl( p, q, -t \bigr) = \frac{ [-t]^{2} }{ 2pq } + \frac{[-t]}{2} \left( \frac{1}{p} + \frac{1}{q} + \frac{1}{pq} \right) + \frac{1}{4} + \frac{1}{12} \left( \frac{p}{q} + \frac{1}{pq} + \frac{q}{p} \right) \\  
                    &\mbox{}& \qquad - \left( \left( \frac{ [-t] }{ p } \right) \right) - \left( \left( \frac{ [-t] }{ q } \right) \right) - \frac{ ((-t)) }{2} \left( \frac 1 p + \frac 1 q \right) - \frac 1 {4p} - \frac 1 {4q} \ . \end{eqnarray*}  
\end{theorem}  
{\it Proof.} We will denote the fractional part of $x$ by $ \{ x \} = x - [x] $. Let $ t = n + r $,  
where $ n \in \Z $ and $ r \in \R , 0 \leq r < 1 $ (so $ r = \{ t \} $), then  
  \begin{eqnarray*} &\mbox{}& \sigma \bigl( a, b, t \bigr) = \sum_{k=0}^{b-1} \left( \left( \frac{ k }{ b } \right) \right) \left( \left\{ \frac{ ka + n + r }{ b } \right\} - \frac 1 2 \right) = \sum_{k=0}^{b-1} \left( \left( \frac{ k }{ b } \right) \right) \left( \left\{ \frac{ ka + n }{ b } \right\} + \frac r b - \frac 1 2 \right) \\  
                    &\mbox{}& \qquad = \sum_{k=0}^{b-1} \left( \left( \frac{ k }{ b } \right) \right) \left( \left( \frac{ ka + n }{ b } \right) \right) + \frac r b \sum_{k=0}^{b-1} \left( \left( \frac{ k }{ b } \right) \right) \stackrel{ (\ref{trivial}) }{ = } \sigma ( a,b,n ) - \frac r {2b} \ . \end{eqnarray*}  
Hence  
  \[ \sigma ( a,b,t ) = \sigma ( a, b, [t] ) - \frac{ \{ t \} }{ 2b } = \sigma ( a, b, [t] ) - \frac 1 {2b} (( t )) - \frac 1 {4b} \ . \]  
Now we can use Corollary \ref{gesvar} for $ \sigma ( a, b, [t] ) $, and the statement follows.  
\hfill {} $\Box$  
 
\vspace{1cm} 
{\bf Acknowledgements} We thank the referees for insightful comments on the first version of this paper.  


\nocite{*}
\addcontentsline{toc}{subsubsection}{References} 
\bibliography{thesis} 
\bibliographystyle{alpha}


 \sc Department of Mathematical Sciences\\ 
 State University of New York\\ 
 Binghamton, NY 13902--6000 \\ 
 {\tt matthias@math.binghamton.edu}\\

 Department of Mathematics\\
 Temple University\\
 Philadelphia, PA 19122 \\
 {\tt srobins@math.temple.edu}

\end{document}